\documentclass[12pt,english,table]{article}
\usepackage{lmodern}

\usepackage[T1]{fontenc}
\usepackage{textcomp}
\usepackage[latin9]{inputenc}
\usepackage{float}
\usepackage{amsmath}
\usepackage{amsthm}
\usepackage{amssymb}
\usepackage{geometry}
\usepackage{setspace}
\usepackage{esint}
\usepackage[authoryear]{natbib}
\makeatletter
\theoremstyle{plain}
\newtheorem{thm}{\protect\theoremname}

\setcitestyle{round}
\usepackage{lscape}
\usepackage{longtable}
\usepackage{xcolor}
\usepackage{rotating}

\usepackage{array}
\usepackage{tikz}
\usepackage{tabularx}\usepackage{multirow}\usepackage{booktabs}
\usepackage{ragged2e}\newcolumntype{C}[1]{>{\centering\arraybackslash}p{#1}}
\usepackage{rotfloat}

\newcolumntype{J}[1]{>{\justify\arraybackslash}p{#1}}
\newcolumntype{R}[1]{>{\RaggedLeft\arraybackslash}p{#1}}
\newcolumntype{Q}[1]{>{\columncolor{Gray}\RaggedLeft\arraybackslash}p{#1}}
\newcolumntype{L}[1]{>{\RaggedRight\arraybackslash}p{#1}}
\newcolumntype{G}{@{\extracolsep{0.5cm}}l@{\extracolsep{0pt}}}%
\newcolumntype{P}[1]{>{\centering\arraybackslash}p{#1}}
\newcolumntype{Y}{>{\centering\arraybackslash}X}
\newcommand{\nhphantom}[1]{\sbox0{#1}\hspace{-\the\wd0}} % command to add negative space, like off-set a minussign

\AtBeginDocument{

}

\usepackage[colorlinks,
            linkcolor=blue,
            anchorcolor=blue,
            citecolor=blue]{hyperref}

\usepackage{babel}

\makeatother

\usepackage{babel}
\providecommand{\theoremname}{Theorem}

\begin{document}
\title{A Unifying Integral Representation of the Gamma Function and Its Reciprocal\thanks{We are grateful to Christian Berg, Steen Thorbjørnsen, and Boris Hanin
for valuable comments.}{\normalsize\emph{\medskip{}
}}}
\author{\textbf{Peter Reinhard Hansen}$^{\mathsection}$$\quad$\textbf{ }and$\quad$\textbf{Chen
Tong}$^{\ddagger}$\thanks{Corresponding author: Chen Tong, Email: tongchen@xmu.edu.cn. Chen
Tong acknowledges financial support from the Youth Fund of the National
Natural Science Foundation of China (72301227), and the Ministry of
Education of China, Humanities and Social Sciences Youth Fund (22YJC790117).}\\
 {\normalsize$^{\mathsection}$}{\normalsize\emph{Department of Economics,
University of North Carolina at Chapel Hill}}\\
{\normalsize$^{\ddagger}$}{\normalsize\emph{Department of Finance,
School of Economics \& Wang Yanan Institute}}\\
{\normalsize\emph{ for Studies in Economics (WISE), Xiamen University\medskip{}
 }}}
\date{{\normalsize\emph{\today}}}
\maketitle
\begin{abstract}
We derive an integral expression $G(z)$ for the reciprocal gamma
function, $1/\Gamma(z)=G(z)/\pi$, that is valid for all $z\in\mathbb{C}$,
without the need for analytic continuation. The same integral avoids
the singularities of the gamma function and satisfies $G(1-z)=\Gamma(z)\sin(\pi z)$
for all $z\in\mathbb{C}$.

\bigskip{}
\bigskip{}
\end{abstract}
{\small\textit{{\noindent}Keywords:}}{\small{} Reciprocal gamma function,
Gamma function, Integral representation, Complex analysis, }Digamma
function, Laplace integral.

\noindent{\small\textit{{\noindent}JEL Classification:}}{\small{}
C02}{\small\par}

\clearpage{}

\section{Introduction}

The gamma function is defined by $\Gamma(z)\equiv\int_{0}^{\infty}t^{z-1}e^{-t}\mathrm{d}t,$
for $\operatorname{Re}(z)>0$, and extended to the entire complex
plane by analytic continuation, satisfying $\Gamma(z+1)=z\Gamma(z)$.
It has singularities (poles) at $z=0,-1,-2,\ldots$.

The reciprocal gamma function, $1/\Gamma(z)$, has no singularities
and is analytic everywhere on $\mathbb{C}$. \citet{Laplace:1785}
showed that 
\begin{equation}
\frac{1}{\Gamma(z)}=\frac{1}{2\pi}\int_{-\infty}^{+\infty}w^{-z}e^{w}\mathrm{d}t,\qquad\operatorname{Re}(z)>0,\label{eq:Laplace}
\end{equation}
where $w=\sigma+it$, with $\sigma>0$ being an arbitrary positive
constant.\footnote{Laplace described (\ref{eq:Laplace}) as a ``\emph{rèsultat assez
remarquable}'' in \citet[p.264]{Laplace:1785} and a ``\emph{rèsultat
remarquable}'' in \citet[p.134]{Laplace:1812}.} The expression can be found in \citet[§63.8]{Nielsen:1906}, where
the requirement, $\operatorname{Re}(z)>0$, is explicitly noted, and
\citet{Pribitkin:2002} demonstrates that a wide range of results
for the gamma function can be derived from this expression. Both \citet{Nielsen:1906}
and \citet{Pribitkin:2002} cite \citet{Laplace:1812} and refer to
the expression as the \emph{Laplace Integral}.\footnote{\citet{Nielsen:1906} cite the 2nd edition of \citet{Laplace:1812},
which is dated (1814).} 

In this paper, we derive a globally valid integral expression for
the reciprocal gamma function. It is similar to (\ref{eq:Laplace}),
but holds for all $z\in\mathbb{C}$. This new integral representation
provides a unified framework for defining both the reciprocal gamma
function and $\Gamma(z)\sin(\pi z)$ without relying on analytic continuation. 

\section{A Globally Valid Integral Expression}

The new expressions for the reciprocal gamma function and gamma function
are the following:
\begin{thm}
\label{Theorem:NewRecipGammaIntegralExp}Let $G(z)\equiv\int_{-\infty}^{\infty}w^{1-2z}e^{w^{2}}\mathrm{d}t$,
where $w=\sigma+it$ with $\sigma>0$, then
\begin{eqnarray}
\frac{1}{\Gamma(z)} & = & \frac{1}{\pi}G(z),\quad\qquad\forall z\in\mathbb{C},\label{eq:HansenTong}\\
\Gamma(z)\sin(\pi z) & = & G(1-z),\qquad\forall z\in\mathbb{C}.\label{eq:gamma2II(1-z)}
\end{eqnarray}
\end{thm}
From the proof of Theorem 1, it follows that the integral, $G(z)$,
is well defined for all $z\in\mathbb{C}$, satisfies the reflection
formula, $G(z)G(1-z)=\pi\sin(\pi z)$ for all $z\in\mathbb{C}$, and
that $\sigma$ is an arbitrary positive constant.

Note that $G(z)$ is more than a simple change of variables, $w\mapsto w^{2}$,
applied to (\ref{eq:Laplace}), because this would also alter the
integration path in $\mathbb{C}$. The key difference between the
integral in Theorem \ref{Theorem:NewRecipGammaIntegralExp} and the
Laplace integral is that $w^{-y}e^{w^{2}/2}\rightarrow0$ as $t\rightarrow\pm\infty$
for all $y\in\mathbb{C}$, whereas $w^{-y}e^{w}$ does not converge
as $t\rightarrow\pm\infty$ for $\operatorname{Re}[y]\leq0$. This
is the reason the new integral representation is valid over the entire
complex plane, while (\ref{eq:Laplace}) only holds for $\operatorname{Re}(z)>0$.\footnote{The integral in Theorem \ref{Theorem:NewRecipGammaIntegralExp} is
not in the comprehensive book by \citet{GradshteynRyzhig:2007}. It
does include (\ref{eq:Laplace}), see \citet[eq. 8.315.2]{GradshteynRyzhig:2007},
albeit their expression contains a typographical error, (``$(a+it)^{2}$''
in the denominator should be ``$(a+it)^{z}$'').}

The integral, $G(z)$, unifies expressions for both the gamma function
and the reciprocal gamma function that are valid across the entire
complex plane without the need for analytic continuation. 

The proof of Theorem 1 utilizes distinct approaches for different
regions of the complex plane. For $\operatorname{Re}(z)>\tfrac{1}{2}$
we can use existing results for moments of a Gaussian random variable,
whereas for $\operatorname{Re}(z)\leq\tfrac{1}{2}$ we will establish
the result by means of a contour integral. 

\noindent \textbf{Proof of Theorem \ref{Theorem:NewRecipGammaIntegralExp}.}
We begin by proving the result for the case $\operatorname{Re}(z)>\tfrac{1}{2}$
using results for a standard Gaussian random variable, $X\sim N(0,1)$.
For $\operatorname{Re}(y)>0$ we have
\begin{eqnarray*}
\mathbb{E}|X|^{y-1} & = & \int_{-\infty}^{\infty}|x|^{y-1}\frac{e^{-x^{2}/2}}{\sqrt{2\pi}}\mathrm{d}x=\sqrt{\tfrac{2}{\pi}}\int_{0}^{\infty}x^{y-1}e^{-x^{2}/2}\mathrm{d}x\\
 & = & \sqrt{\tfrac{2^{y-1}}{\pi}}\int_{0}^{\infty}u^{(y-2)/2}e^{-u}\mathrm{d}u=\tfrac{2^{(y-1)/2}}{\sqrt{\pi}}\Gamma\left(\tfrac{y}{2}\right),
\end{eqnarray*}
where we used the substitution $u=x^{2}/2$. From \citet[theorem 1]{HansenTong:2024MomentsMGF}
we obtain a different expression, $\mathbb{E}|X|^{y-1}=\frac{\Gamma(y)}{\pi}\int_{-\infty}^{+\infty}\ensuremath{\frac{e^{w^{2}/2}}{w^{y}}\mathrm{d}t}=\frac{\Gamma(y)}{\pi}2^{\frac{1-y}{2}}\int_{-\infty}^{+\infty}\ensuremath{\frac{e^{w^{2}/2}}{w^{y}}\mathrm{d}t}$
(also for $\operatorname{Re}(y)>0$). Combining the two expressions
for $\mathbb{E}|X|^{y-1}$, we find
\[
\frac{1}{\pi}\int_{-\infty}^{+\infty}\ensuremath{\frac{e^{w^{2}}}{w^{y}}\mathrm{d}t}=\tfrac{2^{y-1}}{\sqrt{\pi}}\frac{\Gamma\left(\frac{y}{2}\right)}{\Gamma(y)}=\tfrac{2^{y-1}}{\sqrt{\pi}}\frac{\sqrt{\pi}2^{1-y}}{\Gamma\left(\frac{y+1}{2}\right)}=\frac{1}{\Gamma\left(\frac{y+1}{2}\right)},
\]
where we applied the Legendre duplication formula. Setting $y=2z-1$
proves (\ref{eq:HansenTong}) for $\operatorname{Re}(z)>1/2$. 

\begin{figure}[H]
\begin{centering}
\ \hspace{2cm}\begin{tikzpicture}[scale=0.75]
% Define coordinates for points A, B, C, D, E 
\coordinate (A) at (2.83,-2.83);  % Point A
\coordinate (Tneg) at (0,-2.83); % Point -T 
\coordinate (B) at (2.83,2.83); % Point B 
\coordinate (T) at (0,2.83); % Point T 
\coordinate (C) at (0,4);  % Point C 
\coordinate (D) at (0,0);  % Point D (Origin) 
\coordinate (E) at (0,-4);  % Point E

% Axes 
\draw[->, very thick] (-3.5,0) -- (5.5,0) node[right] {$\text{Re}(w)$}; 
\draw[->, very thick] (0,-4.6) -- (0,5) node[above] {$\text{Im}(w)$};
\node[left] at (C) {$R$}; 
	\draw[very thick] (-0.1,4) -- (C);
\node[left] at (0,2.83) {$T$}; 
	\draw[very thick] (-0.1,2.83) -- (0,2.83);
	\draw[dotted] (0,2.83) -- (4,2.83);
\node[left] at (0,-2.83) {$-T$}; 
	\draw[very thick] (-0.1,-2.83) -- (0,-2.83);
	\draw[dotted] (0,-2.83) -- (4,-2.83);
\node[left] at (E) {$-R$}; 
	\draw[very thick] (-0.1,-4) -- (E);
\node at (3.05,-0.23) {$\sigma$}; 
	\draw[very thick] (2.83,0) -- (2.83,-0.1);
	
% five points
\fill[blue] (A) circle (2pt);  % Bullet at A 
    \node[blue] at (3,-3.3) {$A$}; 
	\draw[-, blue] (2.95,-3.1) -- (A);
\fill[blue] (B) circle (2pt);  % Bullet at B 
	\node[blue] at (2.83,3.35) {$B$}; 
	\draw[-, blue] (2.83,3.1) -- (B);
\fill[blue] (C) circle (2pt);  % Bullet at C 
	\node[blue] at (1.0,4.3) {$C$}; 
	\draw[-, blue] (0.7,4.2) -- (C);
\fill[blue] (D) circle (2pt);  % Bullet at D 
	\node[blue] at (-0.5,-0.5) {$D$}; 
	\draw[-,blue] (-0.3,-0.3) -- (D);
\fill[blue] (E) circle (2pt);  % Bullet at E
	\node[blue] at (0.5,-3.5) {$E$}; 
	\draw[-, blue] (0.3,-3.7) -- (E);

% R arrows
\draw[-, dashed] (D) -- (B);
\node[right,black] at (0.9,1.75) {$R$}; 
% Angles
\draw[-,black] (0.65,0.0) arc[start angle=00, end angle=45, radius=0.65] node[pos=0.8, right] {$\Theta$};

% Contour lines
\draw[->,  thick, red] (A) -- (B);
\draw[->,  thick, red] (B) arc[start angle=45, end angle=90, radius=4];
\draw[->, thick, red] (C) -- (D); 
\draw[->, thick, red] (D) -- (E); 
\draw[->, thick, red] (E) arc[start angle=-90, end angle=-45, radius=4];

% Arrows T 
\draw[->|,dashed] (4,1.7) -- (4.0,2.82); 
\draw[->|,dashed] (4,1.1) -- (4,0); 
\draw[->|,dashed] (4,-1.7) -- (4,-2.82); 
\draw[->|,dashed] (4,-1.1) -- (4.0,0.0); 
\node at (4,1.41) {$T$}; 
\node at (4,-1.41) {$T$};

% Arrows R 
%\draw[->|,thick,blue] (-2,2.3) -- (-2.0,4.0); 
%\draw[->|,thick,blue] (-2,1.7) -- (-2.0,0); 
%\draw[->|,thick,blue] (-2,-2.3) -- (-2.0,-4.0); 
%\draw[->|,thick,blue] (-2,-1.7) -- (-2.0,0.0); 
%\node[blue] at (-2,2.0) {$R$}; 
%\node[blue] at (-2,-2.0) {$R$};

\end{tikzpicture} 
\par\end{centering}
\caption{{\small The closed contour, $\mathcal{C}$, for the integral, $\ensuremath{\protect\oint}_{\mathcal{C}}w^{-y}e^{w^{2}}\mathrm{d}w$.\label{fig:contour2}}}
\end{figure}
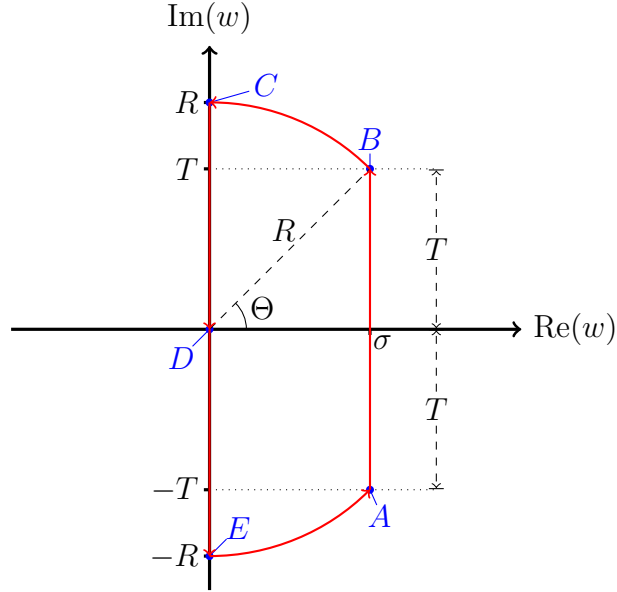

The case $\operatorname{Re}(z)\leq\tfrac{1}{2}$ requires a very different
proof. Define $\tilde{G}(y)\equiv\int w^{-y}e^{w^{2}}\mathrm{d}t=G(\tfrac{y+1}{2})$,
for which (\ref{eq:HansenTong}) and (\ref{eq:gamma2II(1-z)}) can
be restated as,
\[
\tilde{G}(y)=\int_{-\infty}^{+\infty}w^{-y}e^{w^{2}}\mathrm{d}t=\begin{cases}
\cos(\tfrac{\pi}{2}y)\Gamma\left(\tfrac{1-y}{2}\right),\\
\pi\frac{1}{\Gamma\left(\tfrac{1+y}{2}\right)}, & \qquad
\end{cases}\forall y\in\mathbb{C},
\]
and $\tilde{G}(y)\tilde{G}(y+1)=\sqrt{\pi}2^{y}\tilde{G}(2y+1)$,
where we have used Legendre duplication formula. We have $\tilde{G}(y)=\int_{-\infty}^{+\infty}f(y,w)\mathrm{d}t$
where $f(y,w)=w^{-y}e^{w^{2}}$, and 
\[
\int_{-\infty}^{+\infty}f(y,w)\mathrm{d}t=\lim_{T\rightarrow\infty}\int_{-T}^{+T}f(y,w)\mathrm{d}t=\lim_{T\rightarrow\infty}(-i)\int_{\sigma-iT}^{\sigma+iT}f(y,w)\mathrm{d}w.
\]
We proceed to derive and expression for the last integral by means
of a closed contour integral. The path of integration is illustrated
with the red line in Figure \ref{fig:contour2}, where $R=\sqrt{T^{2}+\sigma^{2}}$
and $\Theta=\arccos(\tfrac{\sigma}{R})$ are the polar coordinates
of $(\sigma,T)$.

Because $\operatorname{Re}(y)\leq0$, the function $w\mapsto f(y,w)$
is analytic over the entire region enclosed by the path, such that
$\ensuremath{\oint}_{\mathcal{C}}f(w)\mathrm{d}w=0$ by the Cauchy
Integral Theorem.

Next, we partition the contour integral into the five segments defined
by the points, denoted $A$, $B$, $C$, $D$, and $E$ in Figure
\ref{fig:contour2}. We have
\[
0=\ensuremath{\oint}_{\mathcal{C}}f(w)\mathrm{d}w=I_{AB}(y)+I_{BC}(y)+I_{CD}(y)+I_{DE}(y)+I_{EA}(y),
\]
where the object of interest is the limit of the first integral, $I_{AB}(y)=\int_{\sigma-iT}^{\sigma+iT}\frac{e^{w^{2}}}{w^{y}}\mathrm{d}w$,
as $T\rightarrow\infty$. For the second term, over the $BC$-arc,
we have $w=Re^{i\theta}=R\left(\cos\theta+i\sin\theta\right)$ where
$\theta\in[\Theta,\frac{\pi}{2}]$ and $\mathrm{d}w=iRe^{i\theta}\mathrm{d}\theta$,
such that
\begin{align*}
\left|I_{BC}(y)\right|=\left|\int_{BC}\frac{e^{w^{2}}}{w^{y}}\mathrm{d}w\right| & =\left|\int_{\Theta}^{\pi/2}R^{-y}e^{-i\theta y}e^{R^{2}\cos\left(2\theta\right)+i\left(R^{2}\cos\theta\sin\theta\right)}iRe^{i\theta}d\theta\right|\\
 & \leq\left|R^{1-y}\right|\left|\int_{\Theta}^{\pi/2}e^{R^{2}\cos\left(2\theta\right)}d\theta\right|\\
 & \leq R^{1-{\rm Re}(y)}e^{R^{2}\cos\left(2\Theta\right)}\left|\int_{\Theta}^{\pi/2}d\theta\right|\rightarrow0.
\end{align*}
We used that, $\Theta\rightarrow\pi/2$, $\cos(2\Theta)\rightarrow-1$,
and $R^{p}e^{-R^{2}}\rightarrow0$ for any $p$, as $T\rightarrow\infty$
(because $R\geq T$). Similarly it follows that $\left|I_{EA}(y)\right|\rightarrow0$
as $T\rightarrow\infty$. 

A point on line from $C$ to $D$ is given by $w=e^{\frac{i\pi}{2}}r=ir$
for some $r\in[0,R]$. Hence, 
\begin{eqnarray*}
I_{CD}(y)=\int_{iR}^{0}\frac{e^{w^{2}}}{w^{y}}\mathrm{d}w & = & i\int_{R}^{0}\frac{e^{-r^{2}}}{r^{y}e^{\frac{i\pi}{2}y}}\mathrm{d}r,\\
 & = & -ie^{-\frac{i\pi}{2}y}\int_{0}^{R}\frac{e^{-v}}{v^{\frac{y+1}{2}}}\tfrac{1}{2}\mathrm{d}v,\qquad v=r^{2},\quad\mathrm{d}v=2r\mathrm{d}r\\
 & = & -ie^{-\frac{i\pi}{2}y}\int_{0}^{R}v^{\tfrac{1-y}{2}-1}e^{-v}\tfrac{1}{2}\mathrm{d}v\\
 & \rightarrow & -ie^{-\frac{i\pi}{2}y}\tfrac{1}{2}\Gamma(\tfrac{1-y}{2})\qquad\text{as }R\geq T\rightarrow\infty,
\end{eqnarray*}
where we used $\int_{0}^{\infty}v^{a-1}e^{-v}\mathrm{d}v=\Gamma(a)$,
with $a=\tfrac{1-y}{2}$, which is valid because $\operatorname{Re}(\tfrac{1-y}{2})>0$.
Similarly, $w=-ir$, for $r\in[0,R]$, along the $DE$-line and
\begin{eqnarray*}
I_{DE}(y)=\int_{0}^{-iR}\frac{e^{w^{2}}}{w^{y}}\mathrm{d}w & = & -i\int_{0}^{R}\frac{e^{-r^{2}}}{r^{y}e^{-\frac{i\pi}{2}y}}\mathrm{d}r,\\
 & = & -ie^{\frac{i\pi}{2}y}\int_{0}^{R}\frac{e^{-r^{2}}}{r^{y}}\tfrac{1}{2}\mathrm{d}r\rightarrow-ie^{\frac{i\pi}{2}y}\tfrac{1}{2}\Gamma(\tfrac{1-y}{2}).
\end{eqnarray*}
Thus, for the entire segments from $C$ to $E$, we have
\[
I_{CD}(y)+I_{DE}(y)\rightarrow-i\underset{=\sin(\pi\frac{1-y}{2})}{\underbrace{(\tfrac{1}{2}e^{\frac{i\pi}{2}y}+\tfrac{1}{2}e^{-\frac{i\pi}{2}y})}}\Gamma(\tfrac{1-y}{2}),
\]
and since $\lim_{T\rightarrow\infty}[I_{AB}(y)+I_{CD}(y)+I_{DE}(y)]=0$
we have 
\[
I_{AB}(y)\rightarrow i\sin(\pi\tfrac{1-y}{2})\Gamma(\tfrac{1-y}{2})=i\pi\frac{1}{\Gamma(\tfrac{y+1}{2})},
\]
where we used the Euler reflection formula, $\ensuremath{\Gamma(1-z)\Gamma(z)}=\frac{\pi}{\sin\pi z}$
with $z=\frac{y+1}{2}$, and we have shown that 
\[
\int_{-\infty}^{+\infty}\frac{e^{w^{2}}}{w^{y}}\mathrm{d}t=(-i)\lim_{T\rightarrow\infty}I_{AB}(y)=\pi\frac{1}{\Gamma(\tfrac{y+1}{2})},
\]
for $\operatorname{Re}(y)\leq0$. Since $z=\tfrac{y+1}{2}$ we have
shown (\ref{eq:HansenTong}) for $\operatorname{Re}(z)\leq\tfrac{1}{2}$,
which completes the proof of (\ref{eq:HansenTong}) and (\ref{eq:gamma2II(1-z)})
now follows by Euler's reflection formula.\hfill{}$\square$

\subsection{Digamma function and the Euler-Mascheroni constant}

The derivative of $G(z)$ with respect to $z$ is given by $G^{\prime}(z)\equiv-\int_{-\infty}^{+\infty}w^{1-2z}e^{w^{2}}\log w^{2}\mathrm{d}t$.
Taking the derivative to both sides of the identity, $1/\Gamma(z)=\frac{1}{\pi}G(z)$
and multiplying by $-\Gamma(z)$ leads to the following integral expression
for the digamma function on its domain
\[
\psi(z)=-\frac{G^{\prime}(z)}{G(z)}=\frac{\int_{-\infty}^{+\infty}w^{1-2z}e^{w^{2}}\log w^{2}\mathrm{d}t}{\int_{-\infty}^{+\infty}w^{1-2z}e^{w^{2}}\mathrm{d}t},\qquad z\in\mathbb{C}/\{0,-1,-2,\ldots\}.
\]
As a special case, we have the following expressing for the Euler--Mascheroni
constant, 
\[
\gamma=G^{\prime}(1)/\pi=-\frac{1}{\pi}\int_{-\infty}^{+\infty}\frac{e^{w^{2}}}{w}\log w^{2}\mathrm{d}t=-\frac{1}{\pi}\int_{-\infty}^{+\infty}e^{w^{2}-\frac{1}{2}\log w^{2}}\log w^{2}\mathrm{d}t.
\]

\section{Conclusion}

We have derived a new integral representation for the reciprocal gamma
function that is valid for all $z\in\mathbb{C}$. This representation
complements the Laplace integral, by extending the applicability to
the entire complex plane without requiring analytic continuation.
Furthermore, we showed that the same integral representation provides
an expression for $\Gamma(z)\sin(\pi z)$, unifying the definitions
of the gamma function and its reciprocal with a single integral expression.

The integral, $G(z)$, implies new expressions for other mathematical
objects, as illustrated with the digamma function. In future work
we plan to explore further implications of this result and related
expressions, including connections to other special functions.

\bibliographystyle{apalike}
\bibliography{prh}

\begin{thebibliography}{}

\bibitem[Gradshteyn and Ryzhik, 2007]{GradshteynRyzhig:2007}
Gradshteyn, I.~S. and Ryzhik, I.~M. (2007).
\newblock {\em Table of integrals, series, and products}.
\newblock Academic Press, Amsterdam, 7th edition.

\bibitem[Hansen and Tong, 2024]{HansenTong:2024MomentsMGF}
Hansen, P.~R. and Tong, C. (2024).
\newblock Moments by integrating the moment-generating function.
\newblock {\em arXiv:2410.23587}, [econ.EM].

\bibitem[Laplace, 1785]{Laplace:1785}
Laplace, P.-S. (1785).
\newblock M{\'e}moire sur les approximations des formules qui sont fonctions de
  tr{\`e}s-grands nombres.
\newblock In {\em M{\'e}moires de l'Acad{\'e}mie royale des Sciences de Paris,
  ann{\'e}e 1782}, pages 209--291. Gauthier-Villars, Paris.

\bibitem[Laplace, 1812]{Laplace:1812}
Laplace, P.-S. (1812).
\newblock {\em Th{\'e}orie Analytique des Probabilit\'es}.
\newblock Courcier, Paris.

\bibitem[Nielsen, 1906]{Nielsen:1906}
Nielsen, N. (1906).
\newblock {\em Handbuch Der Theorie Der Gammafunktion}.
\newblock Teubner, Leipzig.

\bibitem[Pribitkin, 2002]{Pribitkin:2002}
Pribitkin, W. d.~A. (2002).
\newblock Laplace's integral, the gamma function, and beyond.
\newblock {\em American Mathematical Monthly}, 109:235--245.

\end{thebibliography}

\end{document}